\documentclass[a4paper,11pt]{amsart}

\usepackage{amsfonts,latexsym,rawfonts,amsmath,amssymb,amsthm, mathrsfs, lscape}
\usepackage{stackrel}
\usepackage[all]{xy}
\usepackage[english]{babel}
\usepackage[utf8]{inputenc}
\usepackage{xfrac}
\usepackage{multirow}
\usepackage{enumitem}
\usepackage[usenames, dvipsnames]{xcolor}
\RequirePackage{tikz, tikz-cd}

\usepackage{array, tabularx}

\usepackage{setspace}
\usepackage{textcomp}

\usepackage[hypertexnames=false,
    pdftex,
    pdfpagemode=UseNone,
    breaklinks=true,
    extension=pdf,
    colorlinks=true,
    linkcolor=blue,
    citecolor=blue,
    urlcolor=blue,
]{hyperref}

\usepackage[top=1.5in, bottom=1in, left=0.9in, right=0.9in]{geometry}

\renewcommand{\Im}{\mathsf{Im}\,}

\newcommand{\ov}{\overline}

\newcommand{\del}{\partial}

\newcommand{\id}{\operatorname{\text{\sf id}}}

\renewcommand{\Im}{\operatorname{Im}}

\newcommand{\LL}{\mathcal{L}}

\newcommand{\al}{\alpha}

\renewcommand{\phi}{\varphi}

\newcommand{\CC}{\mathbb{C}}

\newcommand{\ZZ}{\mathbb{Z}}
\newcommand{\NN}{\mathbb{N}}
\newcommand{\PP}{\mathbb{P}}

\newcommand{\overbar}[1]{\mkern 1.5mu\overline{\mkern-1.5mu#1\mkern-1.5mu}\mkern 1.5mu}

\newcommand*\DA{\mathop{}\!\mathbin\Box}

\def\XXint#1#2#3{{\setbox0=\hbox{$#1{#2#3}{\int}$ }
\vcenter{\hbox{$#2#3$ }}\kern-.6\wd0}}

\allowdisplaybreaks[1]
\newcommand{\introRef}{}

\newtheorem{theorem}{Theorem}[section]
\newtheorem*{theorem*}{Theorem \introRef}
\newtheorem{corollary}[theorem]{Corollary}
\newtheorem*{corollary*}{Corollary \introRef}
\newtheorem{lemma}[theorem]{Lemma}
\newtheorem{definition}[theorem]{Definition}

\newtheorem{remark}[theorem]{Remark}

\newtheorem{example}{Example}[section]

\title{Bott-Chern cohomology of compact Vaisman manifolds}

\author{Nicolina Istrati}
\address{Nicolina Istrati\\FB 12/Mathematik und Informatik\\
Philipps-Universit\"at Marburg\\
Hans-Meer\-wein-Str. 6\\
35032 Marburg\\
Germany}
\email{istrati@mathematik.uni-marburg.de}
\author{Alexandra Otiman}
\address{Alexandra Otiman\\  Institut for Matematik and Aarhus Institute of Advanced Studies, Aarhus University, 8000, Aarhus C, Denmark\\
and\\
Institute of Mathematics “Simion Stoilow” of the Romanian Academy
21 Calea Grivitei Street, 010702, Bucharest, Romania
}
\email{alexandra.otiman@imar.ro, aiotiman@aias.au.dk}
\thanks{The authors are partially supported by a grant of Ministry of Research and Innovation, CNCS - UEFISCDI, project no. PN-III-P1-1.1-TE-2021-0228, within PNCDI III. A. O. is partially supported also by Romanian Ministry of Education and Research, Program PN-III, Project number PN-III-P4-ID-PCE-2020-0025, Contract 30/04.02.2021}
\begin{document}

\begin{abstract}
    We give an explicit description of the Bott-Chern cohomology groups of a compact Vaisman manifold in terms of the basic cohomology. We infer that the Bott-Chern numbers and the Dolbeault numbers of a Vaisman manifold determine each other. On the other hand, we show that the cohomological invariants $\Delta^k$ introduced by Angella-Tomassini are unbounded for Vaisman manifolds. Finally, we give a cohomological characterization of the Dolbeault and Bott-Chern formality for Vaisman metrics.
\end{abstract}
\maketitle

\section{Introduction}

K\"ahler manifolds have remarkable cohomological properties - they are formal, they admit a Hodge decomposition and they satisfy the $\del\ov\del$-lemma. However this is an isolated phenomenon from the point of view of Hermitian geometry, since it is in fact rare that special classes of Hermitian metrics impose good cohomological behaviour. One other such instance is given by the Vaisman manifolds, and the goal of this note is to bring to light some new cohomological behaviour exhibited by this class of manifolds.



Vaisman manifolds, introduced in \cite{V} under the name of generalized Hopf manifolds, are  complex manifolds $X=(M,J)$ admitting a Hermitian metric $\Omega$ satisfying
\[d\Omega=\theta\wedge\Omega\]
where $\theta\neq 0$ is a closed one-form, called the Lee form, which is parallel with respect to the Levi-Civita connection. If we forget the condition that $\theta$ is parallel, then one has the more general class of locally conformally K\"ahler manifolds. 

The vector field $B$ which is metric dual to $\theta$, together with $JB$, define a holomorphic, totally geodesic, transversally K\"ahler foliation $\mathcal F$ with one-dimensional leaves \cite[Theorem 3.1]{V}. This foliation is moreover independent of the choice of the Vaisman metric \cite[Corollary 2.7]{Ts}, and allows one to consider the basic cohomology of the manifold $H_B(\mathcal F)$, which then has the usual K\"ahlerian properties \cite{ka}.

In general, there need not exist relations between the basic cohomology and other standard cohomologies of a manifold endowed with a foliation. However, for Vaisman manifolds this is the case, and to some extent this is the reason why one still has a good control of the cohomological properties of this class of manifolds. 

As an illustration of this fact, the de Rham cohomology of compact Vaisman manifolds has been computed in terms of basic cohomology in \cite{Ka}, see also \cite{V, OV}. Similarly, the Dolbeault cohomology has been computed also in terms of basic cohomology in \cite{T}, see also \cite{Kl}, and in particular it turns out that the Hodge-Fr\"olicher spectral sequence degenerates at the first page, i.e. we have
\[b_k(X)=\sum_{p+q=k}h_{\overline{\partial}}^{p,q}(X).\]

On the other hand, on Vaisman manifolds we do not have the Hodge symmetries $h_{\overline{\partial}}^{p,q}(X)=h_{\overline{\partial}}^{q,p}(X)$ and the $\del\ov\del$-lemma does not hold. This is more generally true on compact locally conformally K\"ahler manifolds \cite{V'}. In particular, the Bott-Chern cohomology 
\[H_{BC}^{\bullet,\bullet}(X):=\frac{\ker\del\cap\ker\ov\del}{\Im\del\ov\del}\]
is then a refined cohomology, which is no longer isomorphic to the de Rham or the Dolbeault cohomology. 

In this note, we compute the Bott-Chern cohomology for any compact Vaisman manifold, in terms of basic cohomology. In order to state our result, let us denote by $\omega=-dJ\theta$ the transverse K\"ahler form associated to a Vaisman metric, let $L:=\omega\wedge\cdot$ and let $\Lambda=L^*$ denote the adjoint operator, viewed as maps acting both on forms and on basic cohomology. We also let $h_0^{p,q}:=\ker\Lambda\cap H^{p,q}_B(\mathcal F)$ denote the primitive basic Hodge numbers of $X$.

\renewcommand{\introRef}{\ref{thm: BC}}
\begin{theorem*}
Let $X$ be an $n$-dimensional compact Vaisman manifold. For any $p,q\in\NN$ we have
\begin{align*}
    H_{BC}^{p,q}(X)\cong H^{p,q}_B(\mathcal F)\cap\ker\Lambda^2&\oplus[\theta^{0,1}]\wedge H_B(\mathcal F)^{p,q-1}\cap\ker L\oplus[\theta^{1,0}]\wedge H_B^{p-1,q}(\mathcal F)\cap\ker L\\
    &\oplus[\theta^{1,0}\wedge\theta^{0,1}]\wedge H_B(\mathcal F)^{p-1,q-1}\cap\ker L.
\end{align*}
In particular, we have the following Bott-Chern numbers for $X$:
\begin{equation*}
    h_{BC}^{p,q}(X)=\begin{cases}h_0^{p,q}+h_0^{p-1,q-1} & p+q<n\\
    h_0^{p-1,q-1}+h_0^{p,q-1}+h_0^{p-1,q} & p+q=n\\
    h_0^{n-p,n-q}+h_0^{n-p-1,n-q}+h_0^{n-p,n-q-1} & p+q>n\end{cases}.
\end{equation*}

\end{theorem*}

This description allows us to deduce that, for complex manifolds of Vaisman type, there is still a close connection between the Dolbeault and the Bott-Chern cohomology. On the one hand, there is a common finite model for both cohomologies:

\renewcommand{\introRef}{\ref{cor: alg}}
\begin{corollary*}
For any compact Vaisman manifold $X$, there exists a finite dimensional bi-graded bi-differential algebra $(A^{\bullet,\bullet},\del_A,\ov\del_A)$ such that 
$$ H_{dR}(X)=H_{dR}(A), \quad H_{\ov\del}(X)=H_{\ov\del}(A), \quad H_{BC}(X)=H_{BC}(A).$$
\end{corollary*}

On the other hand, one cohomology determines the other, at the level of vector spaces:
\renewcommand{\introRef}{\ref{cor: IntCoh}}
\begin{corollary*}
Let $X$ and $Y$ be two compact complex manifolds of Vaisman type. Then
$$ h^{p,q}_{\ov\del}(X)=h^{p,q}_{\ov\del}(Y)\quad \forall p,q\in\NN\quad\Leftrightarrow\quad h^{p,q}_{BC}(X)=h^{p,q}_{BC}(Y) \quad\forall p,q\in\NN.$$
\end{corollary*}

Moreover, while small deformations of Vaisman manifolds don't need to be Vaisman, we can still deduce that the Bott-Chern numbers  of a compact Vaisman manifold are invariant under small analytic deformations (Corollary~\ref{cor: defn}).

These properties exhibit the class of Vaisman manifolds as being cohomologically close to K\"ahler. However, perhaps unexpectedly, we can also notice a behaviour which makes them arbitrarily far from K\"ahler. For each $k\in\NN$, consider the invariants:
\[\Delta^k(X):=\sum_{p+q=k}(h_{BC}^{p,q}(X)+h_{BC}^{n-p,n-q}(X))-2b_k(X).\]
In \cite[Theorems A and B]{AT} it was shown that $\Delta^k(X)\geq 0$ and that a compact complex manifold $X$ satisfies the $\del\ov\del$-lemma precisely when $\Delta^k(X)=0$ for all $k\in\ZZ$. In this sense, the invariants $\Delta^k$ account for how much far $X$ is  from satisfying the $\del\ov\del$-lemma. 

In Corollary \ref{cor: BCI} we compute these invariants for Vaisman manifolds and notice that we have no universal bounds on the $\Delta^k$'s:

\renewcommand{\introRef}{\ref{cor: BCI}}

\begin{corollary*}
Let $n\geq 3$ and $N\in\NN$. Then there is an $n$-dimensional compact Vaisman manifold $X$ with $\Delta^3(X)>N.$
\end{corollary*}

As regards formality of Vaisman manifolds in the sense of Sullivan \cite{sull}, both things can happen: Hopf manifolds are formal, while Kodaira surfaces are not. However, Vaisman manifolds are almost formal, in the sense that their de Rham algebra is an abelian extension of a formal algebra \cite[Theorem 4.8]{ni}.  In this note, we study the metric analogue of formality for Vaisman manifolds, called geometric formality, following \cite{kot}. We recall that a Hermitian metric is called formal, respectively Dolbeault or Bott-Chern formal, if the wedge product of two harmonic forms is again harmonic, where harmonic is meant with respect to the de Rham, respectively Dolbeault or Bott-Chern Laplacian.

Ornea and Pilca  \cite{op} characterized the formality of a Vaisman metric in cohomological terms. Similarly, we do this with respect to the two other  Laplacians. In particular, we observe that while in higher dimension, the formality of a Vaisman metric coincides in the three categories, for complex surfaces we have a new class of manifolds which is Bott-Chern geometrically formal, but not even cohomologically de Rham or Dolbeault formal. 

\renewcommand{\introRef}{\ref{formality}}
\begin{theorem*}
Let $(X, g)$ be a compact $n$-dimensional Vaisman manifold. The following statements are equivalent:
\begin{enumerate}
\item $X$ is  cohomologically Hopf;
    \item $g$ is formal;
    \item $g$ is Dolbeault formal.
    \end{enumerate}

If $n>2$, these are moreover equivalent to 
\begin{enumerate}[label=(4)]
\item $g$ is Bott-Chern formal.
\end{enumerate}
If $n=2$, then there exists a Bott-Chern formal Vaisman metric $g$ on $X$ if and only if $X$ is a diagonal Hopf surface, a Kodaira surface or a finite quotient of these.
 \end{theorem*}

The technical part of our results is Theorem \ref{thm: BC}, for which we give a Hodge theoretical proof. Namely, by \cite[Th\'eor\`eme 2.2]{Sch} there exists a 4th order Laplacian $\Delta_{BC}$ such that 
\[H_{BC}(X)\cong\ker\Delta_{BC}=\ker\del\cap\ker\ov\del\cap\ker(\del\ov\del)^*.\]
We write $\Omega(X)=\Omega(X)_B\otimes\bigwedge\langle\theta^{0,1},\theta^{1,0}\rangle$, where $\Omega(X)_B$ denotes the differentiable forms $\eta$ on $X$ with $\iota_B\eta=\iota_{JB}\eta=0.$ This allows us to translate the equation $\Delta_{BC}\eta=0$ into an equivalent system of equations on each of the $\Omega(X)_B$-components. Using \cite[Theorem 1.1]{Kl}, we deduce moreover that each  $\Omega(X)_B$-component of $\eta$ is in fact a basic form. Using the transverse K\"ahler identities, we are then able to completely characterise the solutions of our system. 

On the way, we also compute the Dolbeault cohomology of Vaisman manifolds in Theorem~\ref{thm: Dolbeault}. This has the advantage, on one hand, to provide a more self contained, and hopefully more transparent, Hodge theoretical proof, which is however similar in spirit to the one of \cite{T}. On the other hand, it illustrates well, but in a much simplified manner, the method of proof of Theorem~\ref{thm: BC}.

\subsection*{Acknowlegements.} We thank Misha Verbitsky for pointing out a previously incomplete argument in the proof of Theorem 6.6, and the anonymous referees for the careful reading and the valuable suggestions.

\section{Preliminaries and Notation}
In this section, we recall some basic notions, some fundamental results and establish some basic facts that we will need in the sequel. Firstly, we present a short account concerning properties of Vaisman manifolds, following \cite{V}. Then, we recall the properties of the transverse K\"ahler structure, following \cite{To, ka}. Finally, we recall and establish some Hodge theoretical facts. In particular, we also fix here the notation used throughout the text. 

\subsection{Vaisman manifolds}
Let $X=(M,J)$ be a compact complex manifold. A Hermitian metric $g$ is called \textit{locally conformally K\"ahler} (LCK) if its fundamental form $\Omega:=g(J\cdot,\cdot)$ satisfies the equation $d\Omega=\theta\wedge\Omega$
for a closed one form $\theta\neq 0$, called the \textit{Lee form} of the metric. The metric dual vector field $B$, defined via $\theta=g(B,\cdot)$ is called the \textit{Lee vector field}.

\begin{remark}
In the present text, we don't allow for $\theta=0$, which would correspond to the definition of K\"ahler metrics. 
\end{remark}

An LCK metric $g$ is called \textit{Vaisman} if
$\nabla^g\theta=0$,
where $\nabla^g$ is the corresponding Levi-Civita connection. It turns out that in this case, the vector fields $B$ and $JB$ are real holomorphic, Killing and have constant norm \cite{V}. Up to multiplying the metric with a positive constant, we can assume that $\theta(B)=||B||=1$, in which case the condition $\LL_B\Omega=0$ also translates into
\begin{equation}\label{Vmetric}
\Omega=-dJ\theta+\theta\wedge J\theta.\end{equation}

In particular, we have $[B,JB]=0$, so the distribution generated by $B$ and $JB$ is integrable and gives rise to a holomorphic foliation $\mathcal F$ with one dimensional parallelizable leaves which are embedded as totally geodesic submanifolds of $X$ \cite[Theorem 3.1]{V}. $\mathcal F$ is usually called \textit{the canonical foliation} of a Vaisman manifold. 

Vaisman manifolds with regular foliation $\mathcal F$ are easy to construct by the following method \cite[Section 5]{V'}. 
\begin{example}\label{ex: Vaisman}
Consider $Y$ a compact projective manifold endowed with a negative holomorphic line bundle $(L,h)$, in the sense that the corresponding Chern curvature $i\Theta(h)$, representing $2\pi c_1(L)$, is a negative form. Let $\lambda\in\CC^*$, $|\lambda|\neq 1$, act on $L$ by fiberwise multiplication. Then the form on $L-0_L$
\[\Omega_s:=\frac{i\del\ov\del||s||^2_h}{||s||^2_h}, \quad s\in L-0_L\]
is the fundamental form of a Vaisman metric, with Lee form $\theta_s=-d\ln||s||^2_h$, which descends to the compact complex manifold $X:=L-0_L/\langle\lambda\rangle$. Moreover, we have a natural holomorphic projection $p:X\rightarrow Y$, whose fibers, isomorphic to $\CC^*/\langle\lambda\rangle$, are precisely the leaves of the foliation $\mathcal F$. In particular, $\mathcal F$ is regular in this case and $X/\mathcal F=Y$.
\end{example}

\subsection{The transverse K\"ahler structure}

Let us recall that given a holomorphic foliation $\mathcal F$ on $X=(M,J)$, a \textit{transverse K\"ahler structure} is given by a Hermitian metric $g_Y$ on the holomorphic bundle $N\mathcal F:=TX/T\mathcal F$ which is holonomy invariant, in the sense that $\LL_V g_Y=0$ for all $V\in\Gamma(T\mathcal F)$, and such that the uniquely associated Chern connection $\nabla^{g_Y}$ of $N\mathcal F$ is torsion free.  The form $\omega_Y:=g_Y(J\cdot,\cdot)$ is a closed basic form, called the \textit{transverse K\"ahler form}. Equivalently, if the foliation $\mathcal F$ is described by the cocycle $(U_{i},p_i,\phi_{ij})$ with local holomorphic transversal $Y$
\begin{equation*}
 \begin{tikzcd}
    U_i\cap U_j\rar{p_i}\dar{p_j} & Y\\
Y \arrow[ru, "\phi_{ij}"']
    \end{tikzcd}
\end{equation*}
where the diagram commutes, $p_i$ are holomorphic submersions and $\phi_{ij}$ are local biholomorphisms, then a transversal K\"ahler structure on $(X,\mathcal F)$ is given by a K\"ahler structure on the local transversal $Y$ which is invariant by the $\phi_{ij}$'s.

In all what follows, $X=(M,J)$ is a compact complex $n$-dimensional manifold endowed with a Vaisman metric $(\Omega,\theta)$ with $||B||=1$. Let us introduce the $(1,1)$-form
\begin{equation}\label{eq: TransvK}
\omega=-dJ\theta=2i\overline{\partial}\theta^{1,0}=-2i\partial\theta^{0,1}.\end{equation}
Then, thanks to eq. \eqref{Vmetric}, $\omega$ is a closed semi-positive form, whose kernel is precisely $T\mathcal F$, and we have:

\begin{theorem}[{\cite[Theorem 3.1]{V}}]
If $(M,J,g,\theta)$ is a Vaisman manifold, then $g$ induces a transverse K\"ahler structure with respect to the canonical foliation $\mathcal F$, with corresponding transverse K\"ahler form $\omega=-dJ\theta$.
\end{theorem}

The basic operators associated to a transversely K\"ahler foliation satisfy the usual K\"ahler relations. First of all, consider the linear bigraded operators 
\begin{gather*}
L:\Omega^{\bullet,\bullet}(X)\rightarrow\Omega^{\bullet+1,\bullet+1}(X),\quad L\alpha:=\omega\wedge\alpha\\ \Lambda:\Omega^{\bullet,\bullet}(X)\rightarrow \Omega^{\bullet-1,\bullet-1}(X), \quad \Lambda=L^*.\end{gather*} 
Then, since $\omega$ has constant rank $n-1$, we have
\begin{equation}\label{LaL}
[\Lambda,L]=(n-1-\deg )\id.
\end{equation}

Next, we define 
\[\Omega(X)_B:=\{\alpha\in\Omega(X)\mid\iota_{B^{1,0}}\alpha=\iota_{B^{0,1}}\alpha=0\}\subset\Omega(X)\]
and the space of basic forms on $X$
\[\Omega_B:=\{\alpha\in\Omega(X)_B\mid \LL_{B^{1,0}}\al=\LL_{B^{0,1}}\alpha=0\}.\]
Then the differential operators $d,\del$ and $\ov\del$ preserve $\Omega_B$, and we set $d_B:=d|_{\Omega_B}$ and similarly for $\del_B$, $\ov\del_B$.

The Vaisman metric induces a natural $L^2$-scalar product on $\Omega(X)$, which we denote by $\langle\cdot,\cdot\rangle$, and we denote by $d^*,\del^*,\ov\del^*, d^*_B, \del_B^*, \ov\del_B^*$ the corresponding formal adjoint operators with respect to this scalar product. Then we have the transverse K\"ahler identities \cite[Section~3.4]{ka}: 

\begin{gather}\label{eq: KahlerIdL}
[L,\overline{\partial}_B^*]=-i\partial_B, \quad [L,\partial_B^*]=i\overline{\partial}_B,\\
\label{eq: KajlerIdLambda}[\Lambda,\overline{\partial}_B]=-i\partial_B^*, \quad [\Lambda,\partial_B]=i\overline\partial_B^*.
\end{gather}

\subsection{Hodge theory}

Recall that, by considering the Laplacians and their corresponding spaces of harmonic forms:
\begin{gather*}
\DA:=[\del,\del^*],\quad 
\overline{\DA\vphantom{A} }:=[\ov\del,\ov\del^*], \quad \Delta:=[d,d^*],\\ 
\mathcal H_\del:=\ker\DA, \quad \mathcal{H}_{\ov\del}:=\ker\overline{\DA\vphantom{A} },\quad \mathcal{H}_\Delta:=\ker\Delta,\quad
\end{gather*}
then, via the standard Hodge theory, we have graded isomorphisms of vector spaces:
\[ H_{\ov\del}(X)\cong \mathcal H_{\ov\del}, \quad H(X,\CC)\cong\mathcal H_\Delta\]
and $L^2$-orthogonal decompositions
\[\Omega(X)=\mathcal H_{\ov\del}\oplus\Im\ov\del\oplus\Im\ov\del^*=\mathcal H_{\del}\oplus\Im\del\oplus\Im\del^*.\]

Similarly, consider the associated basic cohomologies 
\[H^k_B(\mathcal F):=\frac{\ker d_B:\Omega^k_B\rightarrow \Omega^{k+1}_B}{\Im d_B:\Omega^{k-1}_B\rightarrow\Omega^k_B},\quad H_B^{p,q}(\mathcal F):=\frac{\ker \ov\del_B:\Omega^{p,q}_B\rightarrow\Omega^{p,q+1}_B}{\Im\ov\del_B:\Omega^{p,q-1}_B\rightarrow\Omega^{p,q}_B}\]
together with the associated basic Laplacians 
\begin{gather*}
\DA_B:=[\del_B,\del_B^*],\quad \overline{\DA\vphantom{A} }_B:=[\ov\del_B,\ov\del_B^*],\quad \Delta_B:=[d_B,d_B^*].
\end{gather*}
Then, thanks to the transverse K\"ahler identities, one proves in the usual way that 
\[\Delta_B=2\DA_B=2\overline{\DA\vphantom{A} }_B.\] 
We again have \cite[Section 4]{ka'}:
\begin{gather*}
    \mathcal{H}_B:=\ker\overline{\DA\vphantom{A} }_B=\ker\DA_B=\ker\Delta_B, \\
    H_B(\mathcal F)\cong\mathcal H_B, \quad H^k_B(\mathcal F)=\oplus_{p+q=k}H_B^{p,q}(\mathcal F),\\
    \Omega_B=\mathcal H_B\oplus\Im \del_B\oplus\Im\del_B^*=\mathcal H_B\oplus\Im\ov\del_B\oplus\Im\ov\del^*_B.
    \end{gather*}
Moreover, following \cite[Section 3.4.7]{ka}, the operators $L$ and $\Lambda$ act naturally on $H_B(\mathcal F)$ and satisfy the transverse version of the Hard Lefschetz theorem. We introduce, for later use:    
\[ \mathcal{H}_0:=\mathcal{H}_B\cap\ker\Lambda,\quad h_0^{p,q}:=\dim_\CC\mathcal H_0^{p,q}.\]

We have the following immediate fact:

\begin{lemma}\label{lem: Baction}
We have $\LL_{B^{1,0}}=[\del,\iota_{B^{1,0}}]=-\LL_{B^{0,1}}^*$, $\LL_B^*=-\LL_B$ and $\LL_{JB}^*=-\LL_{JB}$. Moreover, $\LL_{B^{1,0}}$ and $\LL_{B^{0,1}}$ act trivially on $\mathcal H_{\ov\del}$, while $\LL_B$ and $\LL_{JB}$ act trivially on $\mathcal H_{\Delta}$.
\end{lemma}
\begin{proof}
Since $B^{1,0}$ is a holomorphic vector field, we have $[\ov\del,\iota_{B^{1,0}}]=0$, thus indeed  $\LL_{B^{1,0}}=[\del,\iota_{B^{1,0}}]$. Since $B^{1,0}$, $B$ and $JB$ are Killing vector fields, we have $\LL^*_{B^{1,0}}=-\ov{\LL_{B^{1,0}}}=-\LL_{B^{0,1}}$, $\LL_B^*=-\LL_B$ and $\LL_{JB}^*=-\LL_{JB}$.

In particular, as $B^{1,0}$ is both Killing and holomorphic, $\LL_{B^{1,0}}$ and its adjoint act on $\mathcal H_{\ov\del}$. Consider $\alpha\in\mathcal H_{\ov\del}$.
Then $\LL_{B^{1,0}}\alpha$, $\LL_{B^{0,1}}\alpha\in\mathcal H_{\ov\del}$. But at the same time we have:
\[\LL_{B^{0,1}}\alpha=\ov\del\iota_{B^{0,1}}\alpha, \quad \LL_{B^{1,0}}\alpha=-\ov\del^*(\theta^{1,0}\wedge\alpha)\]
hence we must have $\LL_{B^{1,0}}\alpha=\LL_{B^{0,1}}\alpha=0$. 

In a similar way we see that $\LL_B$ and $\LL_{JB}$ act trivially on $\mathcal H_\Delta$.


\end{proof}

Recall also the Bott-Chern cohomology, defined as the naturally bigraded algebra
\[ H^{\bullet,\bullet}_{BC}(X):=\frac{\ker\del\cap\ker\ov\del}{\Im\del\ov\del}.\]

Following \cite[Th\'eor\`eme 2.2]{Sch}, there exists a $4$-th order elliptic self-adjoint differential operator 
\[ \Delta_{BC}:=(\partial \overline{\partial})(\partial \overline{\partial})^*+(\partial \overline{\partial})^*(\partial \overline{\partial})+(\overline{\partial}^*\partial)(\overline{\partial}^*\partial)^*+(\overline{\partial}^*\partial)^*(\overline{\partial}^*\partial)+\overline{\partial}^*\overline{\partial}+\partial^*\partial\]
such that
\begin{gather}
H_{BC}(X)\cong\mathcal{H}_{BC}:=\ker\Delta_{BC}=\ker\del\cap\ker\ov\del\cap\ker(\del\ov\del)^*\\
\Omega(X)=\mathcal{H}_{BC}\oplus\Im{\del\ov\del}\oplus(\Im\del^*+\Im\ov\del^*)\label{BC decomp}.
\end{gather}

Finally, let us remark that any form $\alpha\in\Omega(X)$ can be uniquely decomposed as
\begin{equation}\label{theta decomp}
\alpha=\alpha_1+\theta^{0,1}\wedge\alpha_2+\theta^{1,0}\wedge\alpha_3+\theta^{1,0}\wedge\theta^{0,1}\wedge\alpha_4\end{equation}
with $\alpha_i\in\Omega(X)_B$, $i=1,4$, via
\begin{gather*}
    \alpha_4=4\iota_{B^{0,1}}\iota_{B^{1,0}}\alpha\\
    \alpha_3=2\iota_{B^{1,0}}(\alpha-\theta^{1,0}\wedge\theta^{0,1}\wedge\alpha_4)\\
     \alpha_2=2\iota_{B^{0,1}}(\alpha-\theta^{1,0}\wedge\theta^{0,1}\wedge\alpha_4)\\
     \alpha_1=\alpha-\theta^{0,1}\wedge\alpha_2+\theta^{1,0}\wedge\alpha_3+\theta^{1,0}\wedge\theta^{0,1}\wedge\alpha_4.
\end{gather*}
In other words, we have an orthogonal decomposition
\begin{equation*}
\Omega(X)=\Omega(X)_B\oplus\theta^{0,1}\wedge\Omega(X)_B\oplus\theta^{1,0}\wedge\Omega(X)_B\oplus\theta^{1,0}\wedge\theta^{0,1}\wedge\Omega(X)_B
\end{equation*}
and we will write the different operators in block form with respect to this decomposition. In this way, using the fact that $\LL_{B^{1,0}}^*=-\LL_{B^{0,1}}$, we have
\begin{gather}
\nonumber\ov\del=\begin{pmatrix}
    \ov\del_B & 0 & -\frac{i}{2}L & 0\\ 2\LL_{B^{0,1}} & -\ov\del_B & 0 &-\frac{i}{2}L\\
    0 & 0 & -\ov\del_B & 0\\
    0 & 0 & -2\LL_{B^{0,1}} & \ov\del_B
    \end{pmatrix}\quad
    \del=\begin{pmatrix}\del_B & \frac{i}{2}L &0 &0\\
    0 & -\del_B & 0 & 0\\
    2\LL_{B^{1,0}} & 0 & -\del_B & -\frac{i}{2}L\\
    0 & 2\LL_{B^{1,0}} &0 &\del_B
    \end{pmatrix}  \\ \label{op block}
    \\
  \nonumber  \ov\del^*=\begin{pmatrix}\ov\del^*_B & -\LL_{B^{1,0}} & 0 & 0\\ 0 &-\ov\del^*_B & 0 & 0\\ i\Lambda &0 & -\ov\del^*_B & \LL_{B^{1,0}}\\
    0 &i\Lambda &0 &\ov\del^*_B\end{pmatrix}\quad
    \del^*=\begin{pmatrix}\del_B^* &0 &-\LL_{B^{0,1}} &0\\ -i\Lambda &-\del_B^* &0 &-\LL_{B^{0,1}}\\ 0 &0 &-\del_B^* &0\\ 0 &0 &i\Lambda &\del_B^*\end{pmatrix}.
\end{gather}

We notice that the Hodge star operator $*$ preserves the decomposition \eqref{theta decomp} and establishes isomorphisms
\begin{align}
  \nonumber  *&:\Omega^{p,q}(X)_B\rightarrow \theta^{1,0}\wedge\theta^{0,1}\wedge\Omega^{n-q-1,n-p-1}(X)_B\\
 \label{Hodgestar}   *&:\theta^{1,0}\wedge\Omega^{p-1,q}(X)_B\rightarrow\theta^{1,0}\wedge\Omega^{n-q-1,n-p}(X)_B\\
 \nonumber *&:\theta^{0,1}\wedge\Omega^{p,q-1}(X)_B\rightarrow\theta^{0,1}\wedge\Omega^{n-q,n-p-1}(X)_B.
\end{align}

\section{Dolbeault and de Rham Cohomology} 

As a warm up, we start by computing the Dolbeault and the de Rham cohomology of a compact Vaisman manifold. This has already been done in \cite{T} and \cite{Ka, V} from a Hodge theoretical viewpoint, see also \cite{OV, Kl} for different proofs. We provide below yet another Hodge theoretical proof, which in our set up follows quite directly. In what follows, $X$ is a fixed compact complex $n$-dimensional manifold endowed with a Vaisman metric and we use the same notation as in the previous section.

\begin{theorem}\label{thm: Dolbeault}
For any $p,q\in\NN$ we have
\[\mathcal H^{p,q}_{\ov\del}=\mathcal H^{p,q}_0\oplus\theta^{0,1}\wedge\mathcal H^{p,q-1}_0\oplus\theta^{1,0}\wedge \mathcal H_B^{p-1,q}\cap\ker L\oplus\theta^{1,0}\wedge\theta^{0,1}\wedge\mathcal H^{p-1,q-1}_B\cap\ker L. \]
In particular, if $h^{p,q}_0:=\dim\mathcal H_0^{p,q}$ denote the primitive basic Hodge numbers, then we have the following Hodge numbers for $X$:
\begin{equation*}
    h_{\overline{\partial}}^{p,q}(X)=\begin{cases}h_0^{p,q}+h_0^{p,q-1} & p+q<n\\
    h_0^{p,q-1}+h_0^{p-1,q} & p+q=n\\
    h_0^{n-p,n-q}+h_0^{n-p-1,n-q} 
    & p+q>n\end{cases}.
\end{equation*}
\end{theorem}
\begin{proof}
Let $\alpha\in\Omega^{p,q}(X)$ be decomposed as in eq. \eqref{theta decomp}. Then, using Lemma \ref{lem: Baction} and the block formulas \eqref{op block}, we find that $\alpha\in\mathcal H_{\ov\del}^{p,q}=\ker\ov\del\cap\ker\ov\del^*$ if and only if
\begin{gather}
   \nonumber \LL_{B^{1,0}}\alpha_i=\LL_{B^{0,1}}\alpha_i=0, \quad i=1,4,\\
   \nonumber \ov\del_B\alpha_3=\ov\del_B\alpha_4=0=\ov\del_B^*\alpha_1=\ov\del_B^*\alpha_2,\\
   \nonumber \ov\del_B\alpha_1=\frac{i}{2}L\alpha_3,\quad \ov\del_B\alpha_2=-\frac{i}{2}L\alpha_4,\\
  \label{eq: kerDolb}  \ov\del^*_B\alpha_3=i\Lambda\alpha_1,\quad \ov\del_B^*\alpha_4=-i\Lambda\alpha_2.
\end{gather}
In particular, one readily sees that any form in the direct sum 
\[\mathcal H^{p,q}_0\oplus\theta^{0,1}\wedge\mathcal H^{p,q-1}_0\oplus\theta^{1,0}\wedge \mathcal H^{p-1,q}\cap\ker L\oplus\theta^{1,0}\wedge\theta^{0,1}\wedge\mathcal H^{p-1,q-1}_B\cap\ker L\]
satisfies eq. \eqref{eq: kerDolb}, and thus is an element of $\mathcal H_{\ov\del}^{p,q}$. 

Conversely, assume the components of $\alpha$ satisfy eq. \eqref{eq: kerDolb}, so in particular they are all basic forms. With respect to the $L^2$-orthogonal decomposition
\[\Omega_B=\mathcal H_B\oplus\Im\ov\del_B\oplus\Im\ov\del_B^*\]
for a form $\beta\in\Omega_B$, we write $\beta^h$ for its projection onto the $\mathcal H_B$-factor and $\beta^P$ for its projection onto the $\Im P$-factor, where $P\in\{\ov\del_B,\ov\del_B^*\}$.

We infer from eq. \eqref{eq: kerDolb} that we have:
\begin{gather*}
\nonumber\alpha_i=\alpha^h_i+\ov\del_B\gamma_i, \quad \gamma_i\in\Im\ov\del_B^*, \quad i=3,4.\end{gather*}
Moreover, the equations
\begin{equation*}
    \ov\del_B\alpha_1=\frac{i}{2}L\alpha_3^h+\ov\del_B(\frac{i}{2}L\gamma_3),\quad  \ov\del_B\alpha_2=-\frac{i}{2}L\alpha_4^h-\ov\del_B(\frac{i}{2}L\gamma_4)
\end{equation*}
imply
\begin{gather}
    \nonumber L\alpha_3^h=0=L\alpha_4^h\\
\label{eq: Dalpha12} \alpha_1^{\ov\del_B^*}=\frac{i}{2}(L\gamma_3)^{\ov\del_B^*},\quad \alpha_2^{\ov\del_B^*}=-\frac{i}{2}(L\gamma_4)^{\ov\del_B^*}.
\end{gather}

The equation $\ov\del_B^*\ov\del_B\gamma_3=i\Lambda\alpha_1$ immediately gives $\Lambda\alpha^h_1=0$. Moreover, taking the $L^2$-product with $\gamma_3\in\Im\ov\del_B^*$, using the fact that $L\gamma_3$ is $L^2$-orthogonal to $\mathcal H_B$ and that $\alpha_1\in\mathcal H_B\oplus\Im\ov\del_B^*$, we find, via eq. \eqref{eq: Dalpha12}:
\begin{align*}
    ||\ov\del_B\gamma_3||^2_{L^2}=\langle\alpha_1,-iL\gamma_3\rangle=\langle\alpha_1^{\ov\del_B^*},-i(L\gamma_3)^{\ov\del_B^*}\rangle=-2||\alpha_1^{\ov\del_B^*}||^2_{L^2}.
\end{align*}
Thus $\alpha_1=\alpha_1^h$ and $\alpha_3=\alpha_3^h$. 

Similarly, the equation $\ov\del_B^*\ov\del_B\gamma_4=-i\Lambda\alpha_2$ gives $\Lambda\alpha^h_2=0$ and, after taking the $L^2$-product with $\gamma_4$:
\begin{equation*}
    ||\ov\del_B\gamma_4||^2_{L^2}=\langle\alpha_2,iL\gamma_4\rangle=-2||\alpha_2^{\ov\del_B^*}||^2_{L^2}
\end{equation*}
hence $\alpha_2=\alpha_2^h$ and $\alpha_4=\alpha_4^h$. This concludes the proof of the first part.

In order to compute the Hodge numbers, for $p+q<n$ we use the fact that $\ker L\cap\mathcal H_B^k=0$ for $k\leq n-1$. For $p+q>n$, we use the fact that $h_{\overline{\partial}}^{p,q}(X)=h_{\overline{\partial}}^{n-q,n-p}(X)$ together with $h_0^{s,t}=h_0^{t,s}$. For $p+q=n$, we have that $\ker L\cap\mathcal H^{n-1}_B=\mathcal H^{n-1}_0$, since by \eqref{LaL}, $\Lambda L=L\Lambda$ on $\mathcal H_B^{n-1}$ and $L|_{\mathcal H^{n-3}_B}$, $\Lambda|_{\mathcal H^{n+1}_B}$ are injective. Together with $\mathcal{H}^n_0=0$ and $\ker L\cap \mathcal H_B^{n-2}=0$, this gives the desired result. 

\end{proof}

\begin{corollary}\label{cor: deRham}
For any $k\in\NN$ we have
\begin{equation*} 
\mathcal{H}^k_\Delta=\mathcal{H}_0^k\oplus\theta\wedge\mathcal H^{k-1}_0\oplus J\theta\wedge\mathcal H_B^{k-1}\cap\ker L\oplus\theta\wedge J\theta\wedge\mathcal H_B^{k-2}\cap\ker L.\end{equation*}
In particular, the Hodge-Fr\"olicher spectral sequence of $X$
 degenerates at the first page and we have $b_k(X)=\sum_{p+q=k}h_{\overline{\partial}}^{p,q}(X)$.\end{corollary}
 \begin{proof}
Using, via Lemma \ref{lem: Baction}, that \[ d^*(\theta\wedge\eta)=-\LL_B\eta-\theta\wedge d^*\eta, \quad d^*(J\theta\wedge\eta)=-\LL_{JB}\eta-J\theta\wedge d^*\eta \quad \forall \eta\in \Omega^k(X)\] 
it is immediate to check that any form in the above direct sum is in $\ker d\cap\ker d^*=\mathcal H_\Delta$. Moreover, one sees that this direct sum is isomorphic to $\oplus_{p+q=k}\mathcal H_{\ov\del}^{p,q}.$ Thus we obtain
\[\sum_{p+q=k}h_{\overline{\partial}}^{p,q}(X)\leq \dim\mathcal H^k_\Delta=b_k(X)\leq\sum_{p+q=k}h_{\overline{\partial}}^{p,q}(X).\]
Hence we have equality in the above, the Hodge-Frölicher spectral sequence degenerates at the first page and the inclusion of the direct sum into $\mathcal H_\Delta$ is in fact an isomorphism.
 \end{proof}

\section{Bott-Chern Cohomology}

This section is devoted to the computation of the Bott-Chern cohomology of Vaisman manifolds. As before, $X$ is a compact complex $n$-dimensional manifold endowed with a fixed Vaisman structure.

\begin{lemma}\label{lem: BCharmonic}
Let $\alpha\in\Omega(X)$ be decomposed as in formula \eqref{theta decomp}. Then $\alpha\in\mathcal{H}_{BC}$ if and only if the following relations hold:

\begin{gather*}
\LL_{B^{1,0}}\alpha_i=\LL_{B^{0,1}}\alpha_i=0, \quad i=1,4,\\
\ov\del_B\alpha_3=\ov\del_B\alpha_4=0=\del_B\alpha_2=\del_B\alpha_4,\\
\ov\del_B\alpha_1=\frac{i}{2}L\alpha_3, \quad \ov\del_B\alpha_2=-\frac{i}{2}L\alpha_4,\\
\del_B\alpha_1=-\frac{i}{2}L\alpha_2,\quad \del_B\alpha_3=-\frac{i}{2}L\alpha_4,\\
(\del_B\ov\del_B)^*\alpha_1=0,\\
(\del_B\ov\del_B)^*\alpha_2=-i\ov\del_B^*\Lambda\alpha_1,\\
(\del_B\ov\del_B)^*\alpha_3=-i\del_B^*\Lambda\alpha_1,\\
(\del_B\ov\del_B)^*\alpha_4=-\Lambda^2\alpha_1+i\del_B^*\Lambda\alpha_2-i\ov\del_B^*\Lambda\alpha_3.
\end{gather*}
\end{lemma}
\begin{proof}
We use the fact that \[\mathcal{H}_{BC}=\ker\del\cap\ker\ov\del\cap\ker(\del\ov\del)^*\]
together with the block formulas \eqref{op block} for the corresponding operators. This readily implies that the above relations for the $\alpha_i$ produce a form in $\mathcal{H}_{BC}$. 

Conversely, let $\alpha\in\mathcal{H}_{BC}$ and let $G$ be the closure of the complex Lie group generated by $B$ and $JB$. Then $G$ acts on $M$ by holomorphic isometries and according to \cite[Theorem 1.1]{Kl}, it induces a trivial action on the Bott-Chern cohomology, which is $G$-equivariantly  isomorphic to $\mathcal{H}_{BC}$. Moreover, since the decomposition \eqref{theta decomp} is unique and both $\theta^{1,0}$ and $\theta^{0,1}$ are $G$-invariant, we obtain that the forms $\alpha_i$ are indeed basic. The rest of the relations follow directly.
\end{proof}

\begin{theorem}\label{thm: BC}
For any $p,q\in\NN$ 
we have
\[\mathcal{H}^{p,q}_{BC}=\mathcal H^{p,q}_B\cap\ker\Lambda^2\oplus\theta^{0,1}\wedge\mathcal{H}_B^{p,q-1}\cap\ker L\oplus\theta^{1,0}\wedge\mathcal H_B^{p-1,q}\cap\ker L\oplus\theta^{1,0}\wedge\theta^{0,1}\wedge\mathcal H_B^{p-1,q-1}\cap\ker L
.\]
In particular, if $h^{p,q}_0:=\dim\mathcal H_0^{p,q}$ denote the primitive basic Hodge numbers, then we have the following Bott-Chern numbers for $X$:
\begin{equation*}
    h_{BC}^{p,q}(X)=\begin{cases}h_0^{p,q}+h_0^{p-1,q-1} & p+q<n\\
    h_0^{p-1,q-1}+h_0^{p,q-1}+h_0^{p-1,q} & p+q=n\\
    h_0^{n-p,n-q}+h_0^{n-p-1,n-q}+h_0^{n-p,n-q-1} & p+q>n\end{cases}.
\end{equation*}
\end{theorem}
\begin{proof}
Lemma \ref{lem: BCharmonic} implies that the above direct sum is contained in $\mathcal{H}^{p,q}_{BC}$. Conversely, let $\al\in\mathcal{H}^{p,q}_{BC}$ be decomposed as in \eqref{theta decomp}. 
Then Lemma \ref{lem: BCharmonic} implies that each component $\alpha_i$, $i=1,4$ of $\alpha$ is in $\Omega_B$. Since the K\"ahler identities \eqref{eq: KahlerIdL} and \eqref{eq: KajlerIdLambda} give, via the Jacobi identity
\[ [\del_B,\ov\del_B^*]=0=[\ov\del_B,\del_B^*]\]
we have an $L^2$-orthogonal decomposition
\begin{equation}\label{eq: BasicDec}
\Omega_B=\mathcal{H}_B\oplus\Im\del_B\ov\del_B\oplus\Im\del_B\ov\del_B^*\oplus\Im\del^*_B\ov\del_B\oplus\Im\del^*_B\ov\del^*_B.
\end{equation}
As before, for a form $\beta\in\Omega_B$, we will write $\beta^h$ for its projection onto the $\mathcal H_B$-factor and $\beta^P$ for its projection onto the $\Im P$-factor.

We note first that the equations in Lemma \ref{lem: BCharmonic} imply that $\alpha_i^h\in\ker L$ for $i=2,3,4$. Moreover, we can write
\begin{gather}
\nonumber \alpha^{\del_B}_2=\del_B\gamma_2 \quad \gamma_2\in\Im\del_B^*=\ker\del_B^\perp\\
\label{eq: primitives} \alpha^{\ov\del_B}_3=\ov\del_B\gamma_3 \quad \gamma_3\in\Im\ov\del_B^*=\ker\ov\del_B^\perp\\
\nonumber \alpha_4-\alpha_4^h=\del_B\ov\del_B\gamma_4 \quad \gamma_4\in\Im(\del_B\ov\del_B)^*=(\ker\del_B\ov\del_B)^\perp.
\end{gather}

Furthermore, Lemma \ref{lem: BCharmonic} gives $\alpha_1+\frac{i}{2}L\gamma_2\in\ker\del_B$, i.e.
\begin{equation}\label{eq: a1del}
\alpha_1^{\del_B^*}=(-\frac{i}{2}L\gamma_2)^{\del_B^*}.
\end{equation}
Similarly we derive
\begin{gather}
\label{eq: a1ovdel} \alpha_1^{\ov\del_B^*}=(\frac{i}{2}L\gamma_3)^{\ov\del_B^*}\\
\label{eq: gamma432} \gamma_2^{\del_B^*\ov\del_B^*}=-\gamma_3^{\del_B^*\ov\del_B^*}=(\frac{i}{2}L\gamma_4)^{\del_B^*\ov\del_B^*}.
\end{gather}

Let us now study the equation $\ov\del_B^*\del_B^*\alpha_2=-i\ov\del_B^*\Lambda\alpha_1$ from Lemma \ref{lem: BCharmonic}. Using that $\del_B^*$ commutes with $\ov\del_B$ and that $\Lambda$ commutes with $\ov\del_B^*$ and  maps $\mathcal H_B$ to itself, this also reads
\[ (\del_B^*\alpha_2)^{\ov\del_B}=\del_B^*\alpha_2^{\ov\del_B}=-i(\Lambda\alpha_1)^{\ov\del_B}=-i(\Lambda\alpha_1^{\ov\del_B})^{\ov\del_B}.\]
Since $\alpha_1\in\ker(\del_B\ov\del_B)^*$, we have that $\alpha_1^{\del_B\ov\del_B}=0$, hence using \eqref{eq: primitives} and  \eqref{eq: a1del} we find
\begin{align*}
\del_B^*\del_B\gamma_2^{\ov\del_B}&=-i(\Lambda\alpha_1^{\ov\del_B\del_B^*})^{\ov\del_B}=^\eqref{eq: a1del}-i\left(\Lambda(-\frac{i}{2}L\gamma_2)^{\ov\del_B\del_B^*}\right)^{\ov\del_B}\\
&=-\frac{1}{2}\left(\Lambda (L\gamma_2^{\ov\del_B\del_B^*})^{\del_B^*}\right)^{\ov\del_B}=-\frac{1}{2}\left(\Lambda L\gamma_2^{\ov\del_B}\right)^{\del_B^*\ov\del_B}.
\end{align*}
Note that we have again used here that $L$ commutes with $\ov\del_B$ and sends $\Im\del_B$ to itself. Taking the $L^2$-product in the above with $\gamma_2^{\ov\del_B}\in\Im\del_B^*\ov\del_B$, we find:
\begin{equation*}
 ||\del_B\gamma_2^{\ov\del_B}||^2_{L^2}=-\frac{1}{2}||L\gamma_2^{\ov\del_B}||^2_{L^2}.\end{equation*}
Since $\gamma_2\in\Im\del_B^*$, this implies that $\gamma_2^{\ov\del_B}=0$.


Precisely in the same manner, using the equation $\ov\del_B^*\del_B^*\alpha_3=-i\del_B^*\Lambda\alpha_1$ from Lemma \ref{lem: BCharmonic}, \eqref{eq: primitives} and \eqref{eq: a1ovdel}, we find that $\gamma_3^{\del_B}=0$. Therefore, via \eqref{eq: gamma432} we obtain:
\begin{equation}\label{eq: gamma23}
\gamma_2=\gamma_2^{\del^*_B\ov\del_B^*}=-\gamma_3^{\del_B^*\ov\del_B^*}=-\gamma_3=\del^*_B\ov\del_B^*s, \ s\in\Im\del_B\ov\del_B=(\ker\del_B^*\ov\del_B^*)^\perp.
\end{equation}

In particular, via \eqref{eq: a1del} and \eqref{eq: a1ovdel} we have
\begin{equation*}
\al_1^{\del_B^*}=-\frac{i}{2}(L\del_B^*\ov\del_B^*s)^{\del_B^*}, \ \al_1^{\ov\del_B^*}=-\frac{i}{2}(L\del_B^*\ov\del_B^*s)^{\ov\del_B^*}.
\end{equation*}
Since also $\al_1^{\del_B\ov\del_B}=0$, we obtain, using again  the transverse K\"ahler identities:
\begin{align*}
\alpha_1-\alpha_1^h&=-\frac{i}{2}(L\del_B^*\ov\del_B^*s)^{\del^*_B+\ov\del_B^*}\\
&=-\frac{i}{2}(\del_B^*\ov\del_B^*Ls-i\del_B^*\del_Bs+i\ov\del_B\ov\del^*_Bs)^{\del^*_B+\ov\del_B^*}.
\end{align*}

Since $s\in\Im\del_B\ov\del_B$, we deduce:
\begin{equation}\label{eq: alpha1}
\al_1=\al_1^h-\frac{i}{2}\del^*_B\ov\del_B^*Ls.
\end{equation}

Finally we study the last equation from Lemma \ref{lem: BCharmonic}, which via \eqref{eq: primitives}, \eqref{eq: alpha1} and \eqref{eq: gamma23} writes:
\begin{equation*}
\ov\del^*_B\del^*_B\del_B\ov\del_B\gamma_4=-\Lambda^2\alpha_1^h+\frac{i}{2}\del_B^*\ov\del_B^*\Lambda^2Ls-i\ov\del_B^*\del_B^*\Lambda\del_B\del_B^*s-i\ov\del_B^*\del_B^*\Lambda\ov\del_B\ov\del_B^*s.
\end{equation*}
This immediately gives $\Lambda^2\alpha_1^h=0$, as announced. Moreover, taking the $L^2$-product with $\gamma_4$ gives:
\begin{align*}
||\del_B\ov\del_B\gamma_4||_{L^2}^2&=-\langle \del_B^*\ov\del_B^*Ls,\frac{i}{2}L^2\gamma_4\rangle+\langle\ov\del_B^*\del_B^*\del_B\del_B^*s,iL\gamma_4\rangle+\langle\ov\del_B^*\del_B^*\ov\del_B\ov\del_B^*s,iL\gamma_4\rangle\\
&\stackrel[\eqref{eq: gamma23}]{\eqref{eq: gamma432}}{=}-\langle \del_B^*\ov\del_B^*Ls,L\del_B^*\ov\del_B^*s\rangle +\langle\ov\del_B^*\del_B^*\del_B\del_B^*s,2\del_B^*\ov\del_B^*s\rangle+\langle\ov\del_B^*\del_B^*\ov\del_B\ov\del_B^*s,2\del_B^*\ov\del_B^*s\rangle\\
&\stackrel{\eqref{eq: KahlerIdL}}{=}-\langle\del_B^*L\ov\del_B^*s+i\del_B^*\del_Bs,\del_B^*L\ov\del_B^*s-i\ov\del_B\ov\del_B^*s\rangle\\
&\quad\quad -2\langle\del_B^*\del_B\del_B^*\ov\del_B^*s,\del_B^*\ov\del_B^*s\rangle-2\langle\ov\del_B^*\ov\del_B\del_B^*\ov\del_B^*s,\del_B^*\ov\del_B^*s\rangle\\
&=-||\del_B^*L\ov\del_B^*s||^2_{L^2}-2||\del_B\del_B^*\ov\del_B^*s||^2_{L^2}-2||\ov\del_B\del_B^*\ov\del_B^*s||^2_{L^2}.
\end{align*}
Thus $\gamma_4\in\ker\del_B\ov\del_B\cap(\ker\del_B\ov\del_B)^\perp=0$. Via \eqref{eq: gamma432},\eqref{eq: gamma23} and \eqref{eq: a1del} we infer
$\gamma_2=\gamma_3=0$ and $\alpha_i=\alpha_i^h$ for all $1\leq i\leq 4$, which concludes the proof of the first part of the statement. 

In order to deduce the Bott-Chern numbers, we notice first that for $p+q<n$ we have 
\[\mathcal H^{p,q}_{BC}=\mathcal H^{p,q}_B\cap\ker\Lambda^2=\mathcal{H}^{p,q}_0\oplus\omega\wedge\mathcal{H}^{p-1,q-1}_0\]
since $L$ is injective on $\Omega_B^k$ for $k<n-1$. 

For $p+q=n$, we use the fact that $\ker L|_{\mathcal H^{n-2}_B}=0$ and that, by \eqref{LaL}, $\Lambda L=L\Lambda$ on $\mathcal H_B^{n-1}$ and $L|_{\mathcal H^{n-3}_B}$, $\Lambda|_{\mathcal H^{n+1}_B}$ are injective, in order to deduce that $\ker L\cap\mathcal H^{n-1}_B=\mathcal H^{n-1}_0$. Since moreover $\Lambda|_{\mathcal H_B^n}$ is an isomorphism, we find 
\[\mathcal H^{p,q}_{BC}=\Lambda^{-1}\mathcal H^{p-1,q-1}_0\oplus\theta^{0,1}\wedge\mathcal H_0^{p,q-1}\oplus\theta^{1,0}\wedge\mathcal H_0^{p-1,q}\]
which gives us the desired Bott-Chern numbers.

Next, using \eqref{Hodgestar} and the fact that $\Lambda=*^{-1}L*$, we find that $*$ gives an isomorphism
\[\mathcal{H}^{p,q}_{BC}=\mathcal H^{n-q,n-p}_0\oplus\theta^{0,1}\wedge\mathcal{H}_0^{n-q,n-p-1}\oplus\theta^{1,0}\wedge\mathcal H_0^{n-q-1,n-p}\oplus\theta^{1,0}\wedge\theta^{0,1}\wedge\mathcal H_B^{n-q-1,n-p-1}\cap\ker L^2
.\]
In particular, for $p+q>n+1$, the term $\ker L^2$ is trivial and this gives, after conjugation, the desired Bott-Chern numbers. For $p+q=n+1$, we use again that  $\ker L\cap\mathcal H_B^{n-1}=\mathcal H_0^{n-1}$. Since $\mathcal{H}_B=\ker \Lambda\oplus\Im L$, it follows that $\mathcal H_B^{n-3}\cap\ker L^2=L^{-1}\mathcal H_0^{n-1}=0$ and this concludes the proof.
\end{proof}


\section{Relations between cohomologies}
In this section, we deduce a few direct consequences of our computations regarding the interactions between the different cohomologies. In the first two corollaries, we notice that although, for Vaisman manifolds, the Dolbeault cohomology is not isomorphic to the Bott-Chern cohomology, they are still strongly related. More specifically, they admit a common finite-dimensional model and they determine each other. In particular, we also deduce that the Bott-Chern numbers of a Vaisman manifold are invariant under small analytic deformations. On the other hand, in the last corollary, we notice that the difference between the two cohomologies can be arbitrarily high. 

We fix, as before, a compact complex manifold $X$ of Vaisman type. Let $H_B(\mathcal F)$ be the basic cohomology algebra of $X$, let $[\omega]\in H_B^2(\mathcal F)$ be the class of a transverse K\"ahler metric and define the following commutative bi-graded bi-differential algebra $(A^{\bullet,\bullet},\del_A,\ov\del_A)$:
\begin{gather*}
    A=H_B(\mathcal F)\otimes\bigwedge \langle u,\overbar u\rangle, \quad \deg u=(1,0), \quad\deg \overbar u=(0,1)\\
\del_A|_{H_B(\mathcal F)}=\ov\del_A|_{H_B(\mathcal F)}=0,\quad  \del_Au=0=\ov\del_A\ov u,\quad \ov\del_A u=[\omega]=-\del_A\ov u.
\end{gather*}

Then Theorem \ref{thm: Dolbeault}, Corollary \ref{cor: deRham} and Theorem \ref{thm: BC} imply:
\begin{corollary}\label{cor: alg}
$(A^{\bullet,\bullet},\del_A,\ov\del_A)$ is a (finitely generated) model for all three cohomologies of $X$, namely
\begin{gather*}
H_{dR}(X)=H_{dR}(A):=H(A^\bullet,\del_A+\ov\del_A)\\
H_{\ov\del}(X)=H_{\ov\del_A}(A^{\cdot,\bullet})\\
H_{BC}(X)=H_{BC}(A):=\frac{\ker\del_A\cap\ker\ov\del_A}{\Im\del_A\ov\del_A}.
\end{gather*}
\end{corollary}

The following result shows that, just as we can recover the Betti numbers from the Hodge numbers of $X$, we can also recover the Bott-Chern numbers in the same way:
\begin{corollary}\label{cor: IntCoh}
Let $X$ be an $n$-dimensional  compact complex manifold of Vaisman type. Then the Bott-Chern numbers of $X$ can be expressed in terms of the Dolbeault numbers, and vice versa. 
\end{corollary}
\begin{proof}
From Theorem~\ref{thm: Dolbeault} we obtain 
\begin{equation}\label{eq: primToDolb}
    h_0^{p,q}=\sum_{k=0}^q(-1)^k h_{\overline{\partial}}^{p,q-k}\end{equation}
for all $0\leq p+q<n$. Similarly, from Theorem \ref{thm: BC} we obtain
\begin{equation}\label{eq:primToBC}
h_0^{p,q}=\sum_{k=0}^{\min(p,q)}(-1)^kh_{BC}^{p-k,q-k}.
\end{equation}
for all $0\leq p+q<n$. Since the primitive basic numbers determine the other cohomology numbers, this gives the desired result.
\end{proof}

\begin{remark}
We do not know if the above statement holds at the level of the cohomology algebras, i.e. if the algebra structure of $H_{\ov\del}(X)$ determines the algebra structure of $H_{BC}(X)$ or conversely.
\end{remark}

While small deformations of Vaisman manifolds don't need to be Vaisman, we can still deduce the following:

\begin{corollary}\label{cor: defn}
 The Bott-Chern numbers  of a compact Vaisman manifold are invariant under small analytic deformations. 
\end{corollary}

\begin{proof}
 Let $\pi: \mathcal{X} \rightarrow \mathbb{B}$ be a  holomorphic deformation of a Vaisman manifold $X_0:= \pi^{-1}(0)$. Let us put $X_t:= \pi^{-1}(t)$ for $t\in\mathbb{B}$. Using the degeneration of the Hodge-Fr\" olicher spectral sequence of Vaisman manifolds at the first page (Corollary~\ref{cor: deRham}), the invariance of Betti numbers under deformations, the Fr\" olicher inequality and the upper semi-continuity of Hodge numbers (see \cite[Theorem 4.12]{bs76}), we have the following, for $t\in\mathbb{B}$ with $|t|$ small enough:
 \begin{equation*}
     b_k(X_t) \leq \sum_{p+q=k} h^{p, q}_{\overline{\partial}}(X_t) \leq \sum_{p+q=k} h^{p, q}_{\overline{\partial}}(X_0)= b_k(X_0)=b_k(X_t),
 \end{equation*}
 which amounts to $h^{p, q}_{\overline{\partial}}(X_0)=h^{p, q}_{\overline{\partial}}(X_t)$, for any $0 \leq p, q \leq n$.  This means that Hodge numbers are invariant under small analytic deformation of Vaisman type manifolds and applying now Corollary \ref{cor: IntCoh} we get the same for Bott-Chern numbers.
\end{proof}

    

Next, we recall the invariants defined in \cite{AT} for any compact complex manifold $X$: 
\[\Delta^k(X):=\sum_{p+q=k}(h_{BC}^{p,q}(X)+h_{BC}^{n-p,n-q}(X))-2b_k(X), \quad k\in\NN.\]

These invariants measure the failure of the $\del\ov\del$-lemma. When $X$ is a manifold of Vaisman type, it is known \cite[Theorem 2.1]{V'} that it does not satisfy the $\del\ov\del$-lemma. In what follows, we see in fact that there is no universal bound on the $\Delta^k$'s as soon as the dimension of the Vaisman manifold is larger than 2. 

\begin{corollary}\label{cor: BCI}
Let $X$ be a compact complex $n$-dimensional Vaisman manifold, let $b_{k,B}:=\dim H^k_B(\mathcal F)$ and let $b_{0k}:=\dim \mathcal H_0^k=\sum_{p+q=k} h^{p,q}_0$. Then we have:
\begin{gather*}
\Delta^k(X)=\begin{cases}
b_{0k-2} & k<n\\
2b_{0k-2} &k=n\\
b_{02n-k-2} & k>n.
\end{cases}\\
\sum_{0\leq k\leq 2n}\Delta^k(X)=2(b_{n-3,B}+b_{n-2,B}).
\end{gather*}
In particular, $\Delta^0(X)=\Delta^1(X)=\Delta^{n-1}(X)=\Delta^{n}(X)=0$ while $\Delta^2(X)=2$ if $n=2$ and $\Delta^2(X)=1$ if $n>2$. 

Moreover, for $n>2$ we have $\Delta^3(X)\in\{b_{1,B},2b_{1,B}\}$. Thus for fixed dimension $n>2$, the invariants $\Delta^k(X)$ are unbounded for $n$-dimensional compact Vaisman manifolds. 
\end{corollary}
\begin{proof}
Since $h_{\overline{\partial}}^{p,q}(X)=h_{\overline{\partial}}^{n-q,n-p}(X)$, $b_k(X)=b_{2n-k}(X)$ and via Corollary \ref{cor: deRham},  we can also write
\[\Delta^k(X)=\sum_{p+q=k}(h_{BC}^{p,q}(X)-h_{\overline{\partial}}^{p,q}(X)+h^{n-p,n-q}_{BC}(X)-h_{\overline{\partial}}^{n-p,n-q}(X)).\]

By Theorem \ref{thm: Dolbeault} and Theorem \ref{thm: BC} we have
\begin{equation*}
    h_{BC}^{p,q}(X)-h_{\overline{\partial}}^{p,q}(X)=\begin{cases}h_0^{p-1,q-1}-h_0^{p,q-1} & p+q<n\\
    h_0^{p-1,q-1} & p+q=n\\
   h_0^{n-p,n-q-1} & p+q>n\end{cases}
\end{equation*}
which then gives the desired $\Delta^k$'s. We deduce the sum of the $\Delta^k$'s by using the fact that $b_{k,B}=\sum_{j\geq 0}b_{0k-2j}$ for $k<n-1$ by Hard Lefschetz.

Finally, for any $n>2$ and any $g\in\NN$ we can consider a curve $C$ of genus $g$ and take $Y=\PP E$ for some rank $n-1$ holomorphic vector bundle $E$ over $C$. Then $Y$ is an $n-1$-dimensional projective manifold. The construction of Example \ref{ex: Vaisman} produces then an $n$-dimensional Vaisman manifold $X$ with $b_{01}(X)=b_{1,B}(X)=b_1(Y)\geq b_1(C)=2g$, which shows the unboundedness of $\Delta^3$.
\end{proof}

\section{Formality of Vaisman metrics}

The notion of formality concerns the algebraic structure of the various cohomology algebras together with the cup product, and it was introduced and developed by Sullivan \cite{sull}. Specifically, a manifold $M$ is called formal if its algebra of differential forms $(\Omega^\bullet_M,d)$ is quasi-isomorphic, as a commutative differential graded algebra, to $(H^\bullet(M),0)$, meaning that there exists a zig-zag of morphisms between the two algebras inducing isomorphisms in cohomology. While there are obstructions to formality, it is in general not easy to show that a manifold is formal. However, if one knew that the space of harmonic forms with respect to a certain Riemannian metric was closed under the wedge product, then the inclusion of this space into both algebras would immediately produce a quasi-isomorphism and prove formality. 
Motivated by this phenomenon, Kotschick introduced in \cite{kot} the following notion:

\begin{definition}[\cite{kot}] A Riemannian metric $g$ is called formal if the wedge product of two harmonic forms is harmonic. 
\end{definition}

The notion readily extends to other types of cohomologies associated to complex manifolds, namely Dolbeault and Bott-Chern and it is usually referred to as {\it geometric formality}.

\begin{definition}[see \cite{tt}, respectively \cite{at}] A Riemannian metric $g$ on a complex manifold  is called Dolbeault formal (resp. Bott-Chern formal) if the exterior product of two $\Delta_{\overline{\partial}}$ (resp. $\Delta_{BC}$)-harmonic forms is $\Delta_{\overline{\partial}}$ (resp. $\Delta_{BC}$)-harmonic.

\end{definition}

Ornea and Pilca investigated in \cite{op} the formality of Vaisman metrics and linked it with a cohomological property. More exactly, they proved:

\begin{theorem}[{\cite[Theorem 3.2]{op}}]
A Vaisman metric on a compact  complex manifold  is formal if and only if the manifold has the Betti numbers of a Hopf manifold, i.e. $b_1=b_{2n-1}=1$ and $b_{k}=0$, for any $2 \leq k \leq 2n-2$.
\end{theorem}

In the following, we investigate Dolbeault and Bott-Chern formality for Vaisman metrics and obtain cohomological characterizations in the same spirit.  But first we introduce the following definition:

\begin{definition}
A compact complex $n$-dimensional manifold $X$ is called \textit{cohomologically Hopf} if it has the (de Rham) cohomology of a Hopf manifold, i.e. 
\[b_0(X)=b_1(X)=b_{2n-1}(X)=b_{2n}(X)=1\]
and all other Betti numbers vanish.
\end{definition}
\begin{remark}
We notice that in the above definition, for a complex manifold of Vaisman type, we can freely interchange de Rham with Dolbeault, Bott-Chern or primitive. Indeed, since $b_k(X)=\sum_{p+q=k}h_{\overline{\partial}}^{p,q}(X)$, $h_{\overline{\partial}}^{0,1}(X)=h_{\overline{\partial}}^{1,0}(X)+1$ and $h_{\overline{\partial}}^{n-1,n}(X)=h_{\overline{\partial}}^{n,n-1}(X)+1$ by Theorem \ref{thm: Dolbeault}, this  gives that $X$ is de Rham cohomologically Hopf exactly when it is Dolbeault cohomologically Hopf. Moreover, Corollary~\ref{cor: IntCoh} allows us to switch between Dolbeault and Bott-Chern numbers, passing via primitive basic cohomology. 
\end{remark}





\begin{theorem}\label{formality}
Let $(X, g)$ be a compact $n$-dimensional Vaisman manifold. The following statements are equivalent:
\begin{enumerate}
\item $X$ is  cohomologically Hopf;
    \item $g$ is formal;
    \item $g$ is Dolbeault formal.
    \end{enumerate}

If $n>2$, these are moreover equivalent to 
\begin{enumerate}[label=(4)]
\item $g$ is Bott-Chern formal.
\end{enumerate}
If $n=2$, then there exists a Bott-Chern formal Vaisman metric $g$ on $X$ if and only if  $X$ is a diagonal Hopf surface, a Kodaira surface or a finite quotient of these.  
\end{theorem}

\begin{proof}

The equivalence $(1) \Leftrightarrow (2)$ is  \cite[Theorem 3.2]{op}.  

We shall prove now that $(3) \Rightarrow (1)$. Let $\alpha \in \mathcal{H}^{p, q}_0$ such that $p+q<n$. Since $\alpha$ is both $\DA_B$ and $\overline{\DA\vphantom{A} }_B$-harmonic by the transverse K\" ahler identities, $\overline{\alpha}$ shares the same property. By Theorem~\ref{thm: Dolbeault}, this implies that both $\alpha$  and $\overline{\alpha}$ are $\overline{\DA\vphantom{A} }$-harmonic, hence by Dolbeault formality, $\alpha \wedge \overline{\alpha}$ is $\overline{\DA\vphantom{A} }$-harmonic. However, $\alpha \wedge \overline{\alpha}$ is a basic form, therefore by Theorem~\ref{thm: Dolbeault} again, this has to be a primitive form. We deduce from here that $\alpha$ vanishes identically by using the following formula:
\begin{equation}\label{lambda}
  0=\Lambda^{p+q} \alpha \wedge \overline{\alpha} = (-1)^{(p+q)(p+q+1)}(p+q)\langle \alpha, \alpha \rangle. 
\end{equation}
This proves that for $0<p+q<n$, $\mathcal{H}^{p, q}_{\overline{\partial}}$ vanishes identically  unless $p=0$ and $q=1$, in which case $\mathcal{H}^{0, 1}_{\overline{\partial}}$ is generated by $\theta^{0, 1}$. By Serre duality, for $2n>p+q>n$, $\mathcal{H}^{p, q}_{\overline{\partial}}$ vanishes except for $p=n-1$ and $q=n$, in which case $\mathcal{H}^{n-1, n}_{\overline{\partial}}$ is generated by $*\theta^{0, 1}$. Similarly,  $h_{\overline{\partial}}^{p, q}=0$ for $p+q=n$ via Theorem~\ref{thm: Dolbeault},  
hence the manifold is (Dolbeault) cohomologically Hopf.

For the converse implication,  we notice that $(1)$ implies that a basis for the $\overline{\DA\vphantom{A} }$-harmonic forms is given by $\{1, \theta^{0, 1}, *\theta^{0, 1}, \mathrm{vol}_g\}$ and consequently $g$ is by straightforward computation Dolbeault formal. 

Let us now assume that $n>2$. Similarly, $(1)$ implies, via Theorem~\ref{thm: BC}, that a basis for the Bott-Chern harmonic forms is $\{1, \omega, \theta^{1,0}\wedge\omega^{n-1}, \theta^{0,1}\wedge\omega^{n-1}, \theta^{1,0}\wedge \theta^{0,1}\wedge\omega^{n-1}\}$. Consequently, it is Bott-Chern formal by straightforward computation, hence $(1)\Rightarrow(4)$. 

For the converse, we benefit again from the explicit description given by Theorem~\ref{thm: BC}. Let $\alpha \in \mathcal{H}_{B}^{p, q} \cap \mathrm{Ker}\, \Lambda^2$, with $2 \leq p+q<n$. Then both $\alpha$ and $\ov\alpha$ are $\DA_B$ and $\overline{\DA\vphantom{A} }_B$-harmonic, implying that $\alpha$ and $\overline{\alpha}$ are both $\Delta_{BC}$-harmonic. By Bott-Chern formality, $\alpha \wedge \overline{\alpha}$ is also $\Delta_{BC}$-harmonic and since it is basic, by Theorem~\ref{thm: BC} we get $\Lambda^2(\alpha \wedge \overline{\alpha})=0$. By the same formula as in \eqref{lambda} we obtain that $\alpha=0$, which also yields $h_0^{p, q}=0$, for any $0<p+q<n$ and therefore concludes that $(M, J)$ is cohomologically Hopf.

Let us now assume that $n=2$. Then the description of Theorem~\ref{thm: BC} implies that $g$ is formal precisely when the transversely K\"ahler form $\omega$ is formal, i.e. $\forall \alpha,\beta\in \mathcal H_B, \alpha\wedge\beta\in\mathcal H_B$. Let us first note that if $b_1(X)>3$, then $\dim\mathcal H^{1,0}_0\geq 2$, hence there exist two linearly independent holomorphic forms $\alpha_1,\alpha_2\in \mathcal H^{1,0}_0\subset H^{1,0}(X)$. If $g$ is formal, then $g(\alpha_i,\alpha_j)$ is constant for $i,j=1,2$. It follows thus that $\alpha_1$,$\alpha_2$ and $\theta^{1,0}$ are sections of $TX^{1,0}$ which are linearly independent at every point, which contradicts the fact that $\dim_\CC X=2$. Hence, if $g$ is formal, then $b_1(X)\in\{1,3\}$, which by \cite[Theorem 1]{bel} is equivalent to $X$ being an elliptic fibration with multiple fibers over an elliptic curve, or a finite quotient of a diagonal Hopf surface or of a Kodaira surface.

If $b_1(X)=1$, then $X$ is homologically Hopf and we have already seen that this implies $g$ being formal. Assume now that $b_1(X)=3$. If $\pi:X\rightarrow C$ is an elliptic fibration with multiple fibers, then the fibers of $\pi$ are the leaves of the canonical foliation $\mathcal F$. Fix a Vaisman metric $g$ on $X$. Any form $0\neq\alpha\in \mathcal H_B^{1,0}$ is of the form $\alpha=\pi^*\beta$ for some $\beta\in H^{1,0}(C)$. It follows that $\alpha$ has zeroes along the multiple fibers and so cannot have constant norm. This shows that $g$ can never be formal. 

The only other possibility for $b_1(X)=3$ is that $X$ is a primary Kodaira surface. In particular, there exists an elliptic principal bundle $p:X\rightarrow E$ over an elliptic curve $E$. Any flat metric $\omega_E$ on $E$ is formal, and if we choose it with integral cohomology class, then  $-\omega_E$ is the curvature of a negative line bundle $(L,h)$ on $E$, such that the construction of Example~\ref{ex: Vaisman} produces a Vaisman metric $g$ on $X$ with $\omega=p^*\omega_E$. Thus $g$ is formal, which concludes the proof. 
\end{proof}


\begin{thebibliography}{100}

\bibitem[AT13]{AT} D. Angella, A. Tomassini, {\it On the $\del\ov\del$-lemma and Bott-Chern cohomology}, Invent. Math. 192 (2013), no. 1, 71--81.

\bibitem[AT15]{at} D. Angella, A. Tomassini, {\it On Bott-Chern cohomology and formality}, J. Geom. Phys. 93 (2015), 52--61.

\bibitem[BS76]{bs76} C. B\u anic\u a, O. St\u an\u a\c sil\u a, {\it Algebraic methods in the global theory of complex spaces}, Editura Academiei, Bucure\c sti and John Wiley $\&$ Sons, London, New York, Sydney (1976).


\bibitem[B00]{bel} F. Belgun, {\it On the metric structure of non-K\"ahler complex surfaces}, Math. Ann. 317 (2000), 1--40.

\bibitem[CNMY]{ni} B. Cappelletti-Montano, A. De Nicola, J.C. Marrero, I. Yudin, {\it Almost formality of quasi-Sasakian and Vaisman manifolds with applications to nilmanifolds}, Isr. J. Math. 241 (2021), 37--87. 



\bibitem[E90]{ka}  A. El Kacimi-Alaoui,  {\it Op\' erateurs transversalement elliptiques sur un feuilletage riemannien et applications}, Compositio Math. 73 (1990), 57--106.

\bibitem[EH86]{ka'} A. El Kacimi-Alaoui, G. Hector, {\it Decomposition de Hodge basique pour un feuilletage riemannien}, Ann. Inst. Fourier 36 (1986), 207--227.


\bibitem[Kl20]{Kl} N. Klemyatin, {\it Dolbeault cohomology
of compact complex manifolds with an action of a complex Lie group}, J. Geom. Phys. 157 (2020), 103823.





\bibitem[Ka80]{Ka} T. Kashiwada, {\it On V-harmonic forms in compact locally conformal Kähler manifolds with the parallel Lee form}, Kodai Math. J. 3 (1980), 70--82.

\bibitem[Ko01]{kot}  D. Kotschick, {\it On products of harmonic forms}, Duke Math. J. 107 (2001), no. 3, 521--531.



\bibitem[OP11]{op} L. Ornea, M. Pilca, {\it Remarks on the product of harmonic forms}, Pacific Journal of Mathematics, 250 (2011), 353--363

\bibitem[OV22]{OV} L. Ornea, M. Verbitsky, {\it Supersymmetry and Hodge theory on Sasakian and Vaisman manifolds}, Manuscripta Math. (2022), https://doi.org/10.1007/s00229-021-01358-8.

\bibitem[Sch07]{Sch} M. Schweitzer, {\it Autour de la cohomologie de Bott-Chern}, Pr\'epublication de l’Institut Fourier no. 703 (2007), arXiv:0709.3528.

\bibitem[S76]{sull} D. P. Sullivan, {\it Cycles for the dynamical study of foliated manifolds and complex manifolds}, Invent. Math. 36 (1976), no. 1, 225--255.


\bibitem[TT14]{tt} A. Tomassini, S. Torelli, {\it On Dolbeault formality and small deformations}, Internat. J. Mat. 25 (2014), n. 11, 1--9. 

\bibitem[To88]{To} P. Tondeur, {\it Foliations on Riemannian manifolds}, Universitext, Springer-Verlag, 1988.

\bibitem[Ts94]{T}  K. Tsukada, {\it Holomorphic forms and holomorphic vector fields on compact generalized Hopf manifolds}, Compositio Math. 93 (1) (1994), 1--22. 

\bibitem[Ts97]{Ts} K. Tsukada, {\it Holomorphic Maps of Compact Generalized Hopf Manifolds}, Geometriae Dedicata 68 (1997), 61--71.

\bibitem[Va80]{V'} I. Vaisman, {\it On Locally and Globally Conformal K\"ahler Manifolds}, Trans. Amer. Math. Soc. 262 (1980), no. 2, 553--542. 

\bibitem[Va82]{V} I. Vaisman, {\it Generalized Hopf manifolds}, Geom. Dedicata 13 (1982), 231--255.



\end{thebibliography}
\end{document}